\documentclass{nm}
\renewcommand\thisnumber{x}
\renewcommand\thisyear {200x}

\renewcommand\thisvolume{xx}

\usepackage{charter}
\usepackage[charter]{mathdesign}

\usepackage{amssymb}
\usepackage{color}
\usepackage{graphicx}
\usepackage{epstopdf}
\usepackage{amsmath}
\usepackage{verbatim}
\usepackage{amsfonts}
\usepackage{multirow}
\usepackage{float}
\usepackage{forloop}
\usepackage[morefloats=36]{morefloats}%
\input epsf

\long
\def\comment#1{}
\def\ad#1{\begin{aligned}#1\end{aligned}}
\def\b#1{\mathbf{#1}}
\def\a#1{\begin{align*}#1\end{align*}}
\def\an#1{\begin{align}#1\end{align}}

\def\d{\operatorname{div}}

\numberwithin{equation}{section}
\numberwithin{table}{section}
\numberwithin{figure}{section}
\def\boxit#1{\vbox{\hrule height1pt \hbox{\vrule width1pt\kern1pt
#1\kern1pt\vrule width1pt}\hrule height1pt }}
\def\lab#1{\boxit{\small #1}\label{#1}}

\def\meqref#1{\boxit{\small #1}\eqref{#1}}
\def\lab#1{\label{#1}}

\def\meqref#1{\eqref{#1}}

\begin{document}


\markboth{W. Cai and J. Hu and S. Zhang}{Divergence-free  $H(div)$ FEM for
Magnetic Induction}
\title{High Order Hierarchical Divergence-free Constrained Transport $H(div)$ Finite Element Method for \\
Magnetic Induction Equation}


\author[Wei Cai and Jun Hu and Shangyou Zhang]{Wei Cai\affil{1}\comma\corrauth and Jun Hu \affil{2} and Shangyou Zhang \affil{3}}
\address{\affilnum{1}\ Department of Mathematics and Statistics, University of
North Carolina at Charlotte, Charlotte, NC 28223, USA.\\
\affilnum{2}\ School of Mathematical Sciences, Peking University, Beijing
100871, P.R. China. \\
\affilnum{3} \ Department of Mathematics, University of
Delaware, Newark, DE 19716,  USA.}
%
%
\emails{{\tt wcai@uncc.edu} (Wei Cai), {\tt hujun@math.pku.edu.cn } (Jun Hu), {\tt szhang@udel.edu } (Shangyou Zhang)}
%


\begin{abstract}
In this paper, we will use the interior functions of an hierarchical basis for high order $BDM_p$ elements to enforce the divergence-free condition of a magnetic field $B$ approximated by the H(div) $BDM_p$ basis. The resulting constrained finite element method can be used to solve magnetic induction equation in MHD equations. The proposed procedure is based on the fact that the scalar $(p-1)$-th order polynomial space on each element can be decomposed as an orthogonal sum of the subspace defined by the divergence of the interior functions of the $p$-th order $BDM_p$ basis and the constant function. Therefore, the interior functions can be used to remove element-wise all higher order terms except the constant in the divergence error of the finite element solution of $B$-field. The constant terms from each element can be then easily corrected using a first order H(div) basis globally. Numerical results for a 3-D magnetic induction equation show the effectiveness of the proposed method in enforcing divergence-free condition of the magnetic field.
\end{abstract}



\maketitle


\section{Introduction}

Numerical modeling of magneto-hydrodynamic fluids has shown that the observance of the zero
divergence of the magnetic field plays an important role in reproducing the
correct physics in plasmas \cite{brackbill}. Various numerical
techniques have been devised to ensure the computed magnetic field to be
divergence-free \cite{toth}. In the early work of \cite{brackbill} a
projection approach was used to correct the magnetic field to have a zero
divergence. A more natural way to satisfy this constraint is through a class
of the so-called constrained transport (CT) numerical methods based on the
ideas in \cite{evans}. In most CT algorithms for the MHD, the surface
averaged magnetic flux over the surface of a 3-D element is used to represent
the magnetic field so normal continuity of the magnetic field can be assured
while the volume averaged conserved quantities
are used for mass, momentum, and energy variables.

In this paper, we will propose a high order transport finite element method
using a recently developed high order hierarchical basis for the $BDM_p$ element \cite{bdm}
 for the magnetic induction equation in the MHD problems. The divergence condition on the $B$
field is enforced through corrections with interior functions in the basis
set such that the global divergence-free condition will be satisfied.

The paper is organized as follows. In section 2, we will present the hierarchical H(div) basis functions in
various modes (edge, face, and interior). In section 3, we will first characterize the divergence of
the interior basis functions for the hierarchical H(div) basis, then, we will introduce a two-step procedure
to remove non-zero divergence in the finite element solution. Numerical test of the proposed
procedure will be carried out for a 3-D magnetic induction equation in Section 4. Finally, a conclusion is given
in Section 5.

\section{ Basis functions for the tetrahedral element}

\label{sec:constructtet3}
In this section we present hierarchical shape functions proposed in \cite{cai13} for the \ H(div)-conforming
tetrahedral $BDM_p$ element on the canonical reference 3-simplex . The shape functions
are grouped into several categories based upon their geometrical entities on
the reference 3-simplex \cite{ac03}. The basis functions in each category are
constructed so that they are also orthonormal within each category on the reference element.

Any point in the 3-simplex $K^{3}$ is uniquely located in terms of the local
coordinate system $(\xi,\eta,\zeta)$. The vertexes are numbered as
$\mathbf{v}_{0}(0,0,0),$ $\mathbf{v}_{1}(1,0,0),$ $\mathbf{v}_{2}(0,1,0),$
$\mathbf{v}_{3}(0,0,1)$. The barycentric coordinates are given as
\begin{equation}
\lambda_{0}:=1-\xi-\eta-\zeta,\quad\lambda_{1}:=\xi,\quad\lambda_{2}%
:=\eta,\quad\lambda_{3}:=\zeta. \label{eqn:barycentric}%
\end{equation}
The directed tangent on a generic edge $\mathbf{e}_{j}=[j_{1},j_{2}]$ is
defined as
\begin{equation}
\tau^{\mathbf{e}_{j}}:=\tau^{\lbrack j_{1},j_{2}]}=\mathbf{v}_{j_{2}%
}-\mathbf{v}_{j_{1}},\quad j_{1}<j_{2}. \label{eqn:edgedef}%
\end{equation}
The edge is parameterized as
\begin{equation}
\gamma_{\mathbf{e}_{j}}:=\lambda_{j_{2}}-\lambda_{j_{1}},\quad j_{1}<j_{2}.
\label{eqn:edgepara}%
\end{equation}
A generic edge can be uniquely identified with
\begin{equation}
\mathbf{e}_{j}:=[j_{1},j_{2}],\quad j_{1}=0,1,2,\quad j_{1}<j_{2}\leq3,\quad
j=j_{1}+j_{2}+\mathrm{sign}(j_{1}), \label{eqn:edgeidentify}%
\end{equation}
where $\mathrm{sign}(0)=0$. Each face on the 3-simplex can be identified by
the associated three vertexes, and is uniquely defined as
\begin{equation}
\mathbf{f}_{j_{1}}:=[j_{2},j_{3},j_{4}],\quad0\leq\{j_{1},j_{2},j_{3}%
,j_{4}\}\leq3,\quad j_{2}<j_{3}<j_{4}. \label{eqn:facedefine}%
\end{equation}

The standard bases in $\mathbb{R}^{n}$ are noted as $\vec{e}_{i}$,
$i=1,\cdots,n$, and $n=\{2,3\}$.

\subsection{Face functions}

\label{sec:face1}

The face functions are further grouped into two categories: edge-based face
functions and face bubble functions.

\bigskip

\noindent{\underline{Edge-based face functions}: $\ $}

\bigskip These functions are associated with the three edges of a certain face
$\mathbf{f}_{j_{1}}$, and by construction all have non-zero normal components
only on the associated face $\mathbf{f}_{j_{1}}$, \emph{i.e.,}
\begin{equation}
\mathbf{n}^{\mathbf{f}_{j_{k}}}\cdot\Phi_{\mathbf{e}[k_{1},k_{2}]}%
^{\mathbf{f}_{j_{1}},i}=0,\quad j_{k}\neq j_{1}, \label{eqn:edgeface11a}%
\end{equation}
where $\mathbf{n}^{\mathbf{f}_{j_{k}}}$ is the unit outward normal vector to
face $\mathbf{f}_{j_{k}}$.

Using the idea of recursion from \cite{ac03}, independent edge-based face
functions are proposed in \cite{cai13} as follows.

For $p=1$, for each edge we have one face function for this edge%
\begin{equation} \label{1b}
\widetilde{\Phi}_{\mathbf{e}[k_{1},k_{2}]}^{\mathbf{f}_{j_{1}},0}%
=\lambda_{k_{1}}\nabla\lambda_{k_{2}}\times\nabla\lambda_{k_{3}},
\end{equation}
and for $p=2$, one additional new basis function can be constructed as%
\begin{equation}
\widetilde{\Phi}_{\mathbf{e}[k_{1},k_{2}]}^{\mathbf{f}_{j_{1}},1}%
=\lambda_{k_{1}}\lambda_{k_{2}}\nabla\lambda_{k_{3}}\times\nabla\lambda
_{k_{1}},
\end{equation}
which can be shown to satisfy the condition (\ref{eqn:edgeface11a}), and and
for $p\geq3,$ the basis functions are given by
\begin{align}
\widetilde{\Phi}_{\mathbf{e}[k_{1},k_{2}]}^{\mathbf{f}_{j_{1}},i+1}  &
\equiv\ell_{i}(\gamma_{\mathbf{e}_{k}})\widetilde{\Phi}_{\mathbf{e}%
[k_{1},k_{2}]}^{\mathbf{f}_{j_{1}},1}+\ell_{i-1}(\gamma_{\mathbf{e}_{k}%
})\widetilde{\Phi}_{\mathbf{e}[k_{1},k_{2}]}^{\mathbf{f}_{j_{1}}%
,0}\label{eqn:edgeface1var}\\
&  =\ell_{i}(\gamma_{\mathbf{e}_{k}})\left[  \lambda_{k_{1}}\lambda_{k_{2}%
}\nabla\lambda_{k_{3}}\times\nabla\lambda_{k_{1}}\right]  +\ell_{i-1}%
(\gamma_{\mathbf{e}_{k}})\left[  \lambda_{k_{1}}\nabla\lambda_{k_{2}}%
\times\nabla\lambda_{k_{3}}\right],\nonumber
\end{align}
for $\;i=1,\cdots,p-2.$
It can be shown numerically that there are exactly $p$ functions that
are independent and whose normal component is non-zero only on the
associated edge $\mathbf{e}_{k}$.

\bigskip

\noindent{\underline{Face bubble functions}:}

\bigskip

The face bubble functions which belong to each specific group are associated
with a particular face $\mathbf{f}_{j_{1}}$. They vanish on all edges of the
reference 3-simplex $K^{3}$, and the normal components of which vanish on
other three faces, \emph{i.e.,}
\begin{equation}
\mathbf{n}^{\mathbf{f}_{j_{k}}}\cdot\Phi_{m,n}^{\mathbf{f}_{j_{1}}}=0,\quad
j_{k}\neq j_{1}. \label{eqn:bubbleface11}%
\end{equation}
The explicit formula is given as
\begin{equation}
\Phi_{m,n}^{\mathbf{f}_{j_{1}}}=\lambda_{j_{2}}\lambda_{j_{3}}\lambda_{j_{4}%
}L_{m,n}\frac{\nabla\lambda_{j_{3}}\times\nabla\lambda_{j_{4}}}{|\nabla
\lambda_{j_{3}}\times\nabla\lambda_{j_{4}}|},
\end{equation}
\[
L_{m,n}=(1-\lambda_{j_{2}})^{m}(1-\lambda_{j_{2}}-\lambda_{j_{3}})^{n}%
P_{m}^{(2n+3,2)}\left(  \frac{2\lambda_{j_{3}}}{1-\lambda_{j_{2}}}-1\right)
P_{n}^{(0,2)}\left(  \frac{2\lambda_{j_{4}}}{1-\lambda_{j_{2}}-\lambda_{j_{3}%
}}-1\right)
\]

and
\begin{equation}
0\leq\{m,n\},m+n\leq p-3. \label{eqn:facebub4}%
\end{equation}
By construction the face bubble functions share an orthonormal property
on the reference 3-simplex $K^{3}$:
\begin{equation}
<\Phi_{m_{1},n_{1}}^{\mathbf{f}_{j_{1}}},\Phi_{m_{2},n_{2}}^{\mathbf{f}%
_{j_{1}}}>|_{K^{3}}=\delta_{m_{1}m_{2}}\delta_{n_{1}n_{2}},
\label{eqn:facebubbleorthonormal1a}%
\end{equation}
for $0\leq
\{m_{1},m_{2},n_{1},n_{2}\},m_{1}+n_{1},m_{2}+n_{2}\leq p-3.$

\subsection{Interior functions}

Interior functions will have zero components on all four faces while some of
them may still have non-zero tangential components (Edge-based and face-based
interior functions defined below).

\label{sec:interior1} The interior functions are classified into three
categories: edge-based, face-based and bubble interior functions. By
construction the normal component of each interior function vanishes on all
faces of the reference 3-simplex $K^{3}$, \emph{i.e.,}
\begin{equation}
\mathbf{n}^{\mathbf{f}_{j}}\cdot\Phi^{\mathbf{t}}=0,\quad j=\{0,1,2,3\}.
\label{eqn:interfun11}%
\end{equation}

\noindent{\underline{Edge-based interior functions}:}

\bigskip

The tangential component of each edge-based function does not vanish on the
associated only edge $\mathbf{e}_{k}:=[k_{1},k_{2}]$ but vanishes all other
five edges, \emph{i.e.,}
\begin{equation}
\tau^{\mathbf{e}_{j}}\cdot\Phi_{\mathbf{e}[k_{1},k_{2}]}^{\mathbf{t}%
,i}=0,\quad\mathbf{e}_{j}\neq\mathbf{e}_{k}, \label{eqn:edgeinterr11}%
\end{equation}
where $\tau^{\mathbf{e}_{j}}$ is the directed tangent along the edge
$\mathbf{e}_{j}:=[j_{1},j_{2}]$. The shape functions are given as
\begin{equation}
\Phi_{\mathbf{e}[k_{1},k_{2}]}^{\mathbf{t},i}=\ \lambda_{k_{1}}\lambda_{k_{2}%
}\left\{  (1-\lambda_{k_{1}})^{i}P_{i}^{(1,2)}\left(  \frac{2\lambda_{k_{2}}%
}{1-\lambda_{k_{1}}}-1\right)  \right\}  \frac{\tau^{\mathbf{e}_{k}}}%
{|\tau^{\mathbf{e}_{k}}|}, \label{edgeInter}%
\end{equation}
where $i=0,1,\cdots,p-2.$

Again one can prove the orthonormal property of edge-based interior functions:
\begin{equation}
\label{eqn:edge1faceortho}<\Phi_{\mathbf{e}[k_{1},k_{2}]}^{\mathbf{t},m},
\Phi_{\mathbf{e}[k_{1},k_{2}]}^{\mathbf{t},n}>|_{K^{3}} = \delta_{mn},
\quad\{m,n\} = 0,1,\cdots,p-2.
\end{equation}

\bigskip

\noindent{\underline{Face-based interior functions}: }

\bigskip

These functions which are associated with a particular face $\mathbf{f}%
_{j_{1}}$ have non-zero tangential components on their associated face only,
and have no contribution to the tangential components on all other three
faces, \emph{i.e.,}
\begin{equation}
\mathbf{n}^{\mathbf{f}_{j_{k}}}\times\Phi_{m,n}^{\mathbf{t},\mathbf{f}_{j_{1}%
}}=\mathbf{0},\quad j_{k}\neq j_{1}. \label{eqn:inter1face11}%
\end{equation}
Further each face-based interior function vanishes on all the edges of the
3-simplex $K^{3}$, \emph{i.e.,}
\begin{equation}
\tau^{\mathbf{e}_{k}}\cdot\Phi_{m,n}^{\mathbf{t},\mathbf{f}_{j_{1}}}=0.
\label{eqn:face22interr11}%
\end{equation}
The formulas of these functions are given as
\begin{equation}
\Phi_{m,n}^{\mathbf{t},\mathbf{f}_{j_{1}}^{1}}=\lambda_{j_{2}}\lambda_{j_{3}%
}\lambda_{j_{4}}L_{mn}\frac{\tau^{\lbrack j_{2},j_{3}]}}{\left\vert
\tau^{\lbrack j_{2},j_{3}]}\right\vert },\text{ \ \ \ }\Phi_{m,n}%
^{\mathbf{t},\mathbf{f}_{j_{1}}^{2}}=\lambda_{j_{2}}\lambda_{j_{3}}%
\lambda_{j_{4}}L_{mn}\frac{\tau^{\lbrack j_{2},j_{4}]}}{\left\vert
\tau^{\lbrack j_{2},j_{4}]}\right\vert } \label{faceInter}%
\end{equation}%
\[
L_{mn}=(1-\lambda_{j_{2}})^{m}(1-\lambda_{j_{2}}-\lambda_{j_{3}})^{n}%
P_{m}^{(2n+3,2)}\left(  \frac{2\lambda_{j_{3}}}{1-\lambda_{j_{2}}}-1\right)
P_{n}^{(0,2)}\left(  \frac{2\lambda_{j_{4}}}{1-\lambda_{j_{2}}-\lambda_{j_{3}%
}}-1\right)
\]
where  $0\leq\{m,n\},m+n\leq p-3$.
The face-based interior functions enjoy the orthonormal property on the
reference 3-simplex $K^{3}$:%
\begin{equation}
\label{eqn:faceinterorthonormal}<\Phi_{m_{1},n_{1}}^{\mathbf{t},
\mathbf{f}_{j_{1}}^{i}},\Phi_{m_{2},n_{2}}^{\mathbf{t}, \mathbf{f}_{j_{1}}%
^{i}}>|_{K^{3}} = \delta_{m_{1} m_{2}}\delta_{n_{1} n_{2}},
\end{equation}
for $i = \{1,2\}, 0
\le\{m_{1}, m_{2}, n_{1}, n_{2}\}, m_{1} + n_{1}, m_{2} + n_{2} \le p-3.$

\bigskip

\noindent{\underline{Interior bubble functions}:}

\bigskip

The interior bubble functions vanish on the entire boundary $\partial K^{3}$
of the reference 3-simplex $K^{3}$. The formulas of these functions are given
as
\begin{equation}
\Phi_{\ell,m,n}^{\mathbf{t},\vec{e}_{i}}=\lambda_{0}\lambda_{1}\lambda
_{2}\lambda_{3}L_{lmn}\vec{e}_{i},\,i=1,2,3, \label{bubInter}%
\end{equation}
\begin{align*}
L_{lmn}  &  =(1-\lambda_{1})^{m}(1-\lambda_{1}-\lambda_{2})^{n}P_{\ell
}^{(2m+2n+8,2)}\left(  2\lambda_{1}-1\right) \\
&  \cdot P_{m}^{(2n+5,2)}\left(  \frac{2\lambda_{2}}{1-\lambda_{1}}-1\right)
P_{n}^{(2,2)}\left(  \frac{2\lambda_{3}}{1-\lambda_{1}-\lambda_{2}}-1\right)
\end{align*}
where
\[
0\leq\{\ell,m,n\},\ell+m+n\leq p-4.
\]
Again, one can show the orthonormal property of the interior bubble functions

%

\[
<\Phi_{\ell_{1},m_{1},n_{1}}^{\mathbf{t},\vec{e}_{i}},\Phi_{\ell_{2}%
,m_{2},n_{2}}^{\mathbf{t},\vec{e}_{j}}>|_{K^{3}}=\delta_{\ell_{1}\ell_{2}%
}\delta_{m_{1}m_{2}}\delta_{n_{1}n_{2}},
\]
where
\[
0\leq\{\ell_{1},\ell_{2},m_{1},m_{2},n_{1},n_{2}\},\ell_{1}+m_{1}+n_{1}%
,\ell_{2}+m_{2}+n_{2}\leq p-4,\,\{i,j\}=1,2,3.
\]

In Table 1 we summarize the decomposition of the space $\left(  \mathbb{P}%
_{p}(K)\right)  ^{3}$ for the $\mathbf{\mathcal{H}}(\mathbf{div})$-conforming
tetrahedral $BDM_p$ element.

\begin{table}[th]
\begin{center}%
\begin{tabular}
[c]{|c|c|}\hline
Decomposition & Dimension\\\hline
Edge-based face functions & $12p$\\\hline
Face bubble functions & $2(p-2)(p-1)$\\\hline
Edge-based interior functions & $6(p-1)$\\\hline
Face-based interior functions & $4(p-2)(p-1)$\\\hline
Interior bubble functions & $(p-3)(p-2)(p-1)/2$\\\hline
Total & $(p+1)(p+2)(p+3)/2=\dim\left(  P_{p}(K)\right)  ^{3}$\\\hline
\end{tabular}
\caption{Decomposition of the $\left(  \mathbb{P}_{p}(K)\right)  ^{3}$ $BDM_p$ tetrahedral finite element space.}%
\label{tab:deptet}%
\end{center}
\end{table}

\section{Enforcing Divergence-free condition for $B$-field}

The magnetic field B in the MHD equations is assumed to be divergence free, however, the time
evolution from a fully discretized finite element method for the magnetic
equation will render the divergence of $B$ to be non-zero at later time. There are many ways
to remove the non-zero divergence in the magnetic field such as the projection
method through Helmholtz decomposition. In this paper, we will use the
interior functions in the H(div) basis set to correct the non-zero divergence
element by element. Due to the vanishing property of the normal components of
the interior basis functions, such a local correction will still keep the
corrected finite element solution in H(div) globally. The ability
of using only the interior functions to reduce the non-divergence error in the
magnetic field is based on the following result.

First let us denote the subspace spanned all the interior functions defined in
(\ref{edgeInter}), (\ref{faceInter}), and (\ref{bubInter}) as%
\begin{equation}
\Sigma_{\rm int}=\operatorname{span}\{\Phi_{i}\}_{i=1}^{n_{i}},\text{ \ \ }n_{i}=\frac{1}%
{2}(p-1)(p+1)(p+2). \label{BasisInt}%
\end{equation}

\begin{lemma} \label{lemma1} $\d \Sigma_{\rm int}=P_{p-1}(K)\backslash
\{1\}=\{1\}^{\perp}.$
\end{lemma}

\begin{proof}
We will prove the result by subspace inclusion argument. First, we will show
$\d \Sigma_{\rm int}\subset\{1\}^{\perp}$. Take any function $\Phi\in\Sigma_{\rm int}%
$, we have%
\begin{equation}
\int_{K}\d \Phi\cdot1dx=\int_{\partial K}\Phi\cdot\mathbf{n}ds=0,\text{ \ i.e.
\ } \d \Phi\perp1, \label{pf1}%
\end{equation}
due to the fact the the normal component of interior function vanishes on the
edge of the element K, thus $\d \Sigma_{\rm int}\subset\{1\}^{\perp}.$

On the other hand, we will show $\{1\}^{\perp}\subset \d \Sigma_{\rm int}$ by
showing that $(\d \Sigma_{\rm int})^{\perp}\subset\{1\}$, instead.
Let $v\in P_{p-1}(K)$ and $v\in(\d \Sigma_{\rm int})^{\perp}$, then we have%
\begin{equation}
\int_{K}v\cdot \d \Phi dx=0,\text{ \ \ for all }\Phi\in\Sigma_{\rm int},
\label{pf2}%
\end{equation}
which gives with an integration by parts
\begin{equation}
\int_{\partial K}\nabla v\cdot\Phi dx=0, \label{pf3}%
\end{equation}
where the vector field $\nabla v\in\left(  P_{p-2}(K)\right)  ^{3}.$

Now take the tangential vectors of the three edges sharing the common vertex
$\mathbf{v}_{0}$,$\tau^{\lbrack0,1]}$,\newline $\tau^{\lbrack0,2]}$,$\tau^{\lbrack0,3]},$ we
can see easily that the following vector functions are also interior functions
(with zero normal components on all faces),%
\begin{equation}
\lambda_{0}\lambda_{1}g_{1}\tau^{\lbrack0,1]},\lambda_{0}\lambda_{2}g_{2}%
\tau^{\lbrack0,2]},\lambda_{0}\lambda_{3}g_{3}\tau^{\lbrack0,3]}\in
\Sigma_{\rm int} \label{pf4}%
\end{equation}
where the scalar functions $\ g_{i}\ ,i=1,2,3$ are polynomials of degree
$(p-2).$ Next, we construct three bi-orthogonal vectors $s_{j}$ with respect
to $\tau^{\lbrack0,i]},i=1,2,3$ with the following property%
\begin{equation}
s_{j}\cdot\tau^{\lbrack0,i]}=\delta_{ij}. \label{pf5}%
\end{equation}

We can express the vector field $\nabla v$ using the basis vector
$s_{j}$ as follows%
\begin{equation}
\nabla v=f_{1}\mathbf{s}_{1}+f_{2}\mathbf{s}_{2}+f_{3}\mathbf{s}_{3},\text{
\ \ }f_{i}\in P_{p-2}(K),i=1,2,3, \label{pf6}%
\end{equation}
which will be substituted into (\ref{pf3}), resulting in%
\begin{equation}
\int_{\partial K}\left(  f_{1}\mathbf{s}_{1}+f_{2}\mathbf{s}_{2}%
+f_{3}\mathbf{s}_{3}\right)  \cdot\Phi dx=0. \label{pf7}%
\end{equation}
By setting $\Phi=\lambda_{0}\lambda_{1}f_{1}\tau^{\lbrack0,1]}$ in the above
identity and using the bi-orthogonality property of (\ref{pf5}), we have%
\begin{equation}
\int_{\partial K}\lambda_{0}\lambda_{1}f_{1}^{2}dx=0, \label{pf8}%
\end{equation}
which implies that $f_{1}=0$ as both $\lambda_{0},\lambda_{1}$ are positive
inside $\overline{K}$. Similar arguement will show that $f_{2}=f_{3}=0$ also hold.

Therefore, we have $\nabla v=0,$ namely, $v=$const., thus $(\d \Sigma_{\rm int})^{\perp}\subset\{1\}.$
\end{proof}

\bigskip

{\noindent \bf Algorithm:} we propose a two-step algorithm to remove the non-divergence in the
numerical solution for the magnetic field $B$.

\begin{itemize}
\item \textbf{Step 1 (local correction)} Element-wise removal of high order
terms in $\d B$.

Due to Lemma \ref{lemma1}, we can use the interior function in $\Sigma_{\rm int}$ to remove
higher order terms in $\d B$. The remaining component in $\d B$ will be a
constant on each element. We proceed to finding a vector function%
\begin{equation}
\Phi=
{\displaystyle\sum\limits_{i=1}^{n_{i}}}
\alpha_{i}\Phi_{i}\in\Sigma_{\rm int}%
\end{equation}
such that%
\begin{equation}
B_{1}=B+\Phi
\end{equation}
\begin{equation}
div(B+\Phi)=c,
\end{equation}
namely, $div(B+\Phi)\in(\d \Sigma_{\rm int})^{\perp}=\{P_{p-1}(K)\setminus
P_{0}\}^{\perp}$, which gives the following linear system for the unique
    $\Phi\in  \Sigma_{\rm int}$%
\an{ \label{local1}%
  \left\{\ad{ (w,\operatorname{div}\Phi)+(w,\operatorname{div}B)& =0
       \quad\forall w\in P_{p-1}(K)\setminus P_{0}, \\
      (\Phi,\Psi)+(v,\operatorname{div}\Psi)&=0 \quad\forall\Psi\in
      \Sigma_{\mathrm{\rm int}}, }\right. }
    where $v\in P_{p-1}(K)\setminus P_{0}$.
Here \eqref{local1} is a local mixed finite element approximation to the
  following    Poisson equation,
 \a{ -\d \operatorname{grad} v &= -\d B \qquad\hbox{in } K, \\
     \frac{\partial v}{\partial \b n} &=0 \qquad \qquad\hbox{on } \partial K. }

\item \textbf{Step 2 (global correction)} remove the constant term in
$\d B_{1}$ in the whole domain.

Due to the result of lemma 1, we will have a residual constant term left in
the corrected magnetic field $B_{1}$, which can only be removed by a global
correction with the first order H(div) basis defined in Section 2. We proceed
as follows by finding a second function $\Phi\in H(div,\Omega)$ using the
first order $H(div)$ basis functions defined in \eqref{1b},
\begin{equation}
\Phi=
{\displaystyle\sum\limits_{i=1}^{N_{1}}}
\beta_{i}\Phi_{i}%
\end{equation}
such that%
\begin{equation}
B^{\ast}=B_{1}+\Phi
\end{equation}%
\begin{equation}
\int_{\Omega}divB^{\ast}\cdot \d \Phi_{j}dx=0,\text{ \ \ for \ }1\leq j\leq
N_{1},
\end{equation}
resulting into the following linear system for $\Phi$, in a similar argument
for equations (\ref{local1}),
\an{\label{global1} \left\{ \ad{
(\Phi,\Psi)+(v,\operatorname{div}\Psi)  & =0\quad\forall\Psi\in
\hbox{span}\{\Phi_{i},\ i=1,...,N_{1}\},\\
(w,\operatorname{div}\Phi)+(w,\operatorname{div}B_{1})  & =0\quad\forall
w\in\{P_{0}(K)\}. } \right.  }

\end{itemize}

\section{Numerical results}

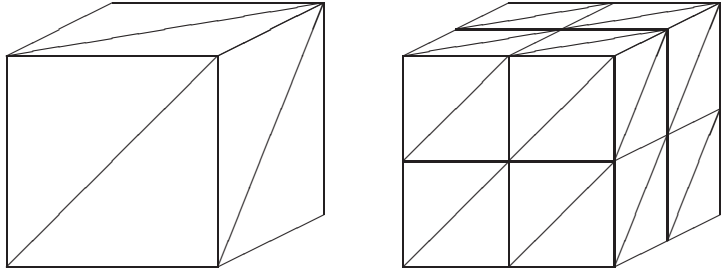
\begin{figure}[tbh]
\begin{center}
\begin{picture}(280,100)(0,0)
\put(0,0){\begin{picture}(100,100)(0,0)
\multiput(0,0)(80,0){2}{\line(0,1){80}}     \multiput(0,0)(0,80){2}{\line(1,0){80}}
\multiput(80,0)(0,80){2}{\line(2,1){40}}    \multiput(0,80)(80,0){2}{\line(2,1){40}}
\multiput(80,0)(40,20){2}{\line(0,1){80}}    \multiput(0,80)(40,20){2}{\line(1,0){80}}
\multiput(0,0)(40,0){1}{\line(1,1){80}}
\multiput(80,0)(0,80){1}{\line(2,5){40}}    \multiput(0,80)(80,0){1}{\line(6,1){120}}
\end{picture} }
\put(150,0){\begin{picture}(100,100)(0,0)
\multiput(0,0)(40,0){3}{\line(0,1){80}}     \multiput(0,0)(0,40){3}{\line(1,0){80}}
\multiput(80,0)(0,40){3}{\line(2,1){20}}    \multiput(0,80)(40,0){3}{\line(2,1){20}}
\multiput(100,10)(0,40){3}{\line(2,1){20}}    \multiput(20,90)(40,0){3}{\line(2,1){20}}
\multiput(80,0)(20,10){3}{\line(0,1){80}}    \multiput(0,80)(20,10){3}{\line(1,0){80}}
\multiput(0,0)(40,0){1}{\line(1,1){80}}     \multiput(40,0)(-40,40){2}{\line(1,1){40}}
\multiput(80,0)(0,80){1}{\line(2,5){40}}    \multiput(0,80)(80,0){1}{\line(6,1){120}}
\multiput(100,10)(-20,30){2}{\line(2,5){20}} \multiput(40,80)(-20,10){2}{\line(6,1){60}}
\end{picture} }
\end{picture}
\end{center}
\caption{The level 1 and 2 uniform grids. }%
\label{criss}%
\end{figure}

We will solve a magnetic induction field equation on the unit cube $\Omega
=[0,1]^{3}$,
\begin{equation}
B_{t}=-\operatorname{div}(BU^{T}-UB^{T}),\label{mhd}%
\end{equation}
with a periodic $B\cdot\mathbf{n}$ boundary condition is considered
\a{ B_1(t,0,y,z)&= B_1(t,1,y,z), \\ B_2(t,x,0,z) &= B_2(t,x,1,z), \\
    B_3(t,x,y,0) &= B_2(t,x,y,1).}
In  (\ref{mhd}),
\begin{equation}
U=%
\begin{pmatrix}
1\\
1\\
0
\end{pmatrix}
.
\end{equation}
Also for  (\ref{mhd}), the initial condition $B(0, \b x )$ is
given by the exact solution
\begin{equation} \label{sol1}
B( t, \b x)=%
\begin{pmatrix}
\sin(2\pi(x+y-z-2t))+\sin(2\pi(y-t))\\
\sin(2\pi(x-t))\\
\sin(2\pi(x+y-z-2t))
\end{pmatrix}
.
\end{equation}
We note that due to the initial condition $\d B( 0,\b x)=0$, \eqref{mhd} ensures
\begin{equation}
\operatorname{div}B(t, \b x )=0,\quad\hbox{ for all }t>0.
\end{equation}

\begin{table}[htb]
\begin{center}
 \begin{tabular}[c]{|c|cc|cc|r|}  
\hline & $ \|B-B_h\|_{0}$ &$h^n$ &
    $ \|\d(B-B_h)\|_0=\|\d B_h \|_0$ & $h^n$  & $\dim V_3$   \\ \hline
2&    0.86157111&0.0&     2.52192914&0.0& 1920\\
3&    0.31631008&1.4&     1.27803825&1.0& 15360\\
4&    0.05191410&2.6&     0.44154994&1.5&122880\\
5&    0.00635164&3.0&     0.07211612&2.6&983040\\
      \hline
  \multicolumn{6}{|c|} {  For $\d B_h=0$ corrected solution with both \meqref{local1} and \meqref{global1}.  }  \\
      \hline
 2&    0.95693378&0.0&     0.00000032&0.0&  1920\\
 3&    0.32810096&1.5&     0.00000012&1.4& 15360\\
 4&    0.05713843&2.5&     0.00000040&0.0& 122880\\
 5&    0.00698271&3.0&     0.00000063&0.0& 983040\\
      \hline
\end{tabular}
 \caption{\lab{b-2} The errors
     and the order of convergence, by the $P_3$ element, for \meqref{mhd} and divergence-free corrections.}
\end{center}

\end{table}

We discritize equation (\ref{mhd}) by a $P_{3}$ $H(\operatorname{div})$ mixed
finite element in Section 2, {\it i.e.} the space $V_{3}$, on the uniform tetrahedral grids depicted in
Fig. \ref{criss}.   For time evolution, we use the
characteristic method. We compute the solution $B(1/2)$ with 100 time
steps with $dt=0.005$. The $L^{2}$ errors are listed in Table
\ref{b-2}, along with the order of convergence, which is $3$ as expected for
the $P_{3}$ mixed finite element. We can see that the numerical solutions are
not divergence-free any more as the divergence-free condition for
$B_{h}(t,\b x)$ is not enforced at the time discretization level.
Following the correction method proposed in Section 3, we correct the solution
$B_{h}(1/2,\b x)$ by direct solution of both \meqref{local1} and \meqref{global1}.
The error and the convergence order are listed also in Table \ref{b-2} at the bottom.
We can see that the divergence of corrected solution is almost zero, up to the computer accuracy.
Moreover, the correction does not change the order of convergence of the finite element
  solution, though the L2 norm of the error for the corrected solution is slightly larger.

\begin{figure}[tbh]
\setlength\unitlength{1in}
\begin{center}
\begin{picture}(3.8,5)(.1,.1)
 \put(0,3.2){ \epsfysize=1.5in \epsfbox{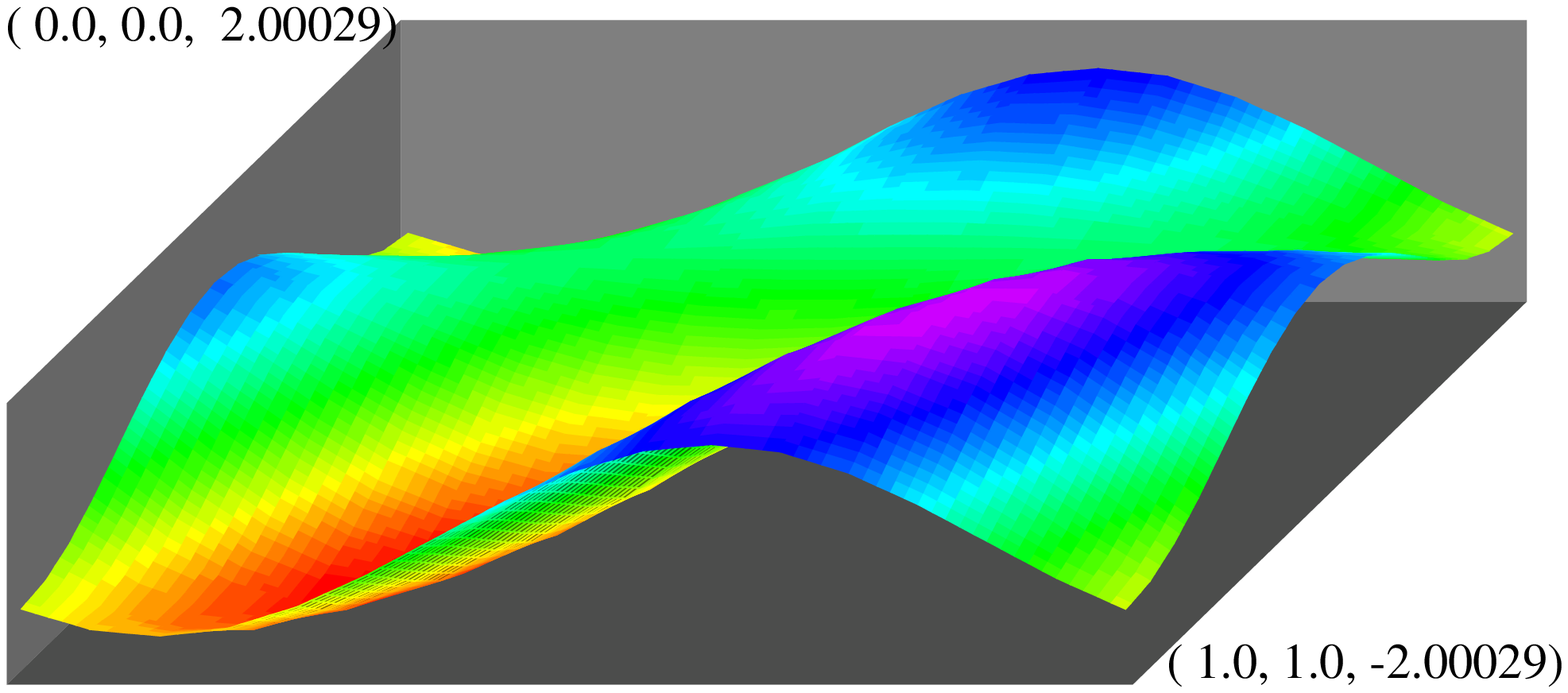}}
  \put(0,1.6){ \epsfysize=1.5in \epsfbox{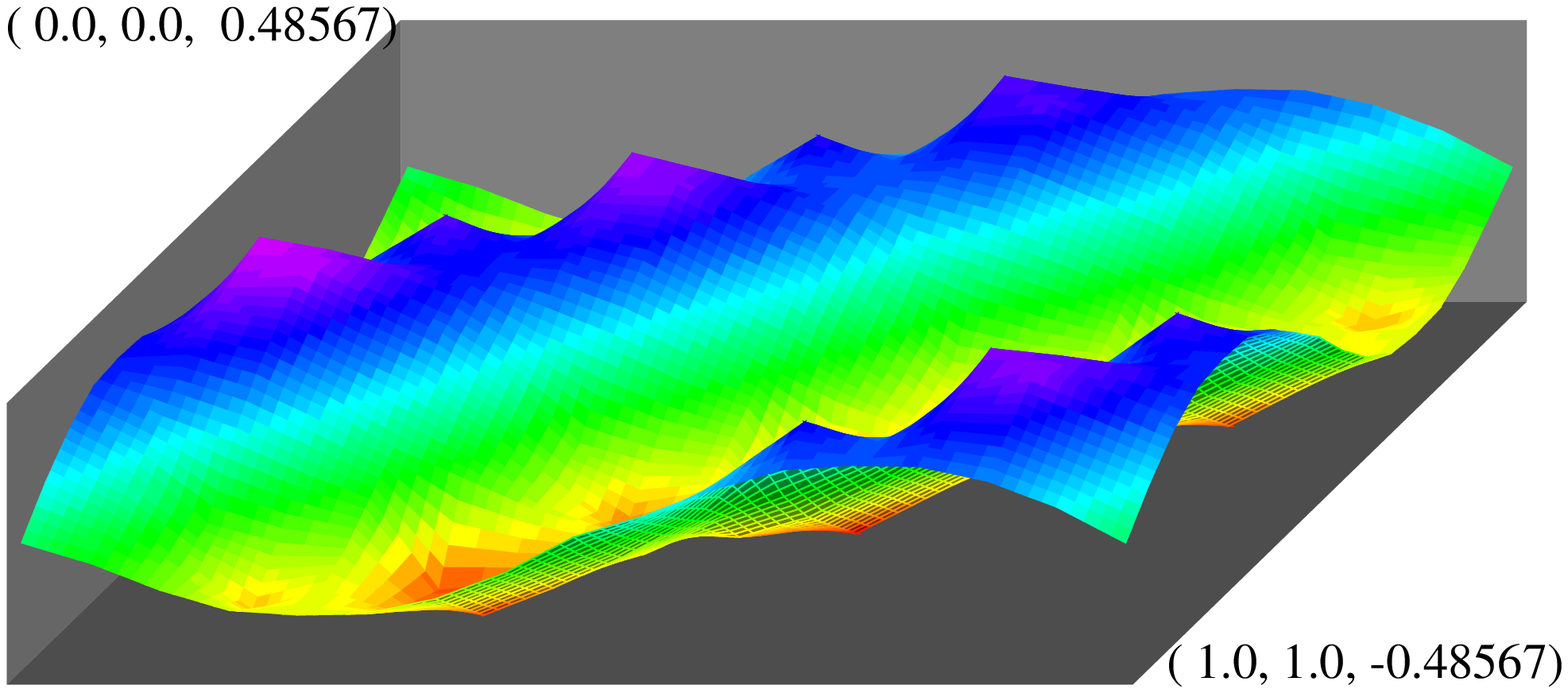}}
 \put(0,0){ \epsfysize=1.5in \epsfbox{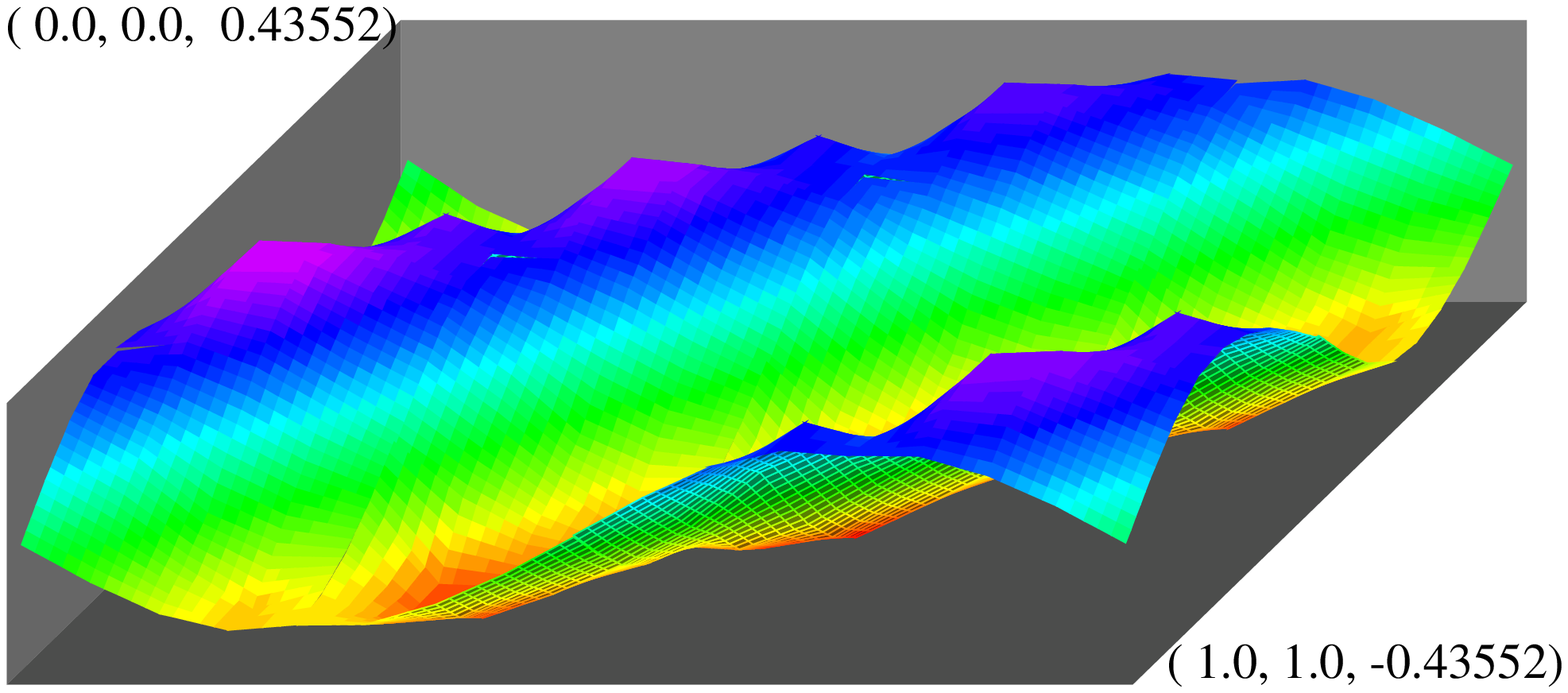}}
\end{picture}
\end{center}
\caption{The solution $ B_1(1/2)$ (top), the error before div-free correction (middle),
   and the error of corrected solution (bottom). }%
\label{f-2}%
\end{figure}

To see the effect of the two correction steps on the solution,
  we plot in Figure \ref{f-2} the solution $ B_1(1/2)$ at $t=1/2$ and the
   errors before and after the correction on a plane at $z=0.485.$
We also compare the new method with the existing, global correction method with the full
$P_3-H(div)$ basis.
To make a comparison to the existing full divergence-zero correction method,
  we list in Table \ref{b-3} (top part) the  errors of the
  divergence-corrected solution with $P_3-H(div)$ basis in \meqref{global1}. The correction equations are
  solved by an iterative Uzawa method, stopped when the $\operatorname{div}$
  norm is less than $0.000005$.
As expected, we see that the computation complexity is much higher than our new, two-step
   correction method.
We notice that from Table \ref{b-3} that the
  corrected solution now has a smaller error while achieving an
  almost zero divergence.
In doing the two-step correction \meqref{local1} and \meqref{global1},
  the work for the global correction \meqref{global1} is of several orders higher than the local correction step.
So, we might wish do just the first step correction \meqref{local1}.
However, the corrected solution is not divergence-free anymore and we list in the bottom half of Table \ref{b-3} the errors and the order of convergence
   for such one-step corrected solutions.

\begin{table}[htb]
\begin{center}
 \begin{tabular}[c]{|c|cc|cc|r|}  
\hline & $ \|B-B_h\|_{0}$ &$h^n$ &
    $ \|\d(B-B_h)\|_0=\|\d B_h \|_0$ & $h^n$  & $\dim V_3$   \\ \hline
  \multicolumn{6}{|c|} {  For a global $\d B_h=0$ corrected solution with  $P_3 - H(div)$ basis.  }  \\
      \hline
 2&    0.84092712&0.0&    0.00000484&0.0&      1920\\
 3&    0.31432784&1.4&    0.00000476&0.0&     15360\\
 4&    0.05155935&2.6&    0.00000487&0.0&    122880\\
 5&    0.00618254&3.0&    0.00000447&0.0&    983040\\
      \hline
  \multicolumn{6}{|c|} {  For local $\d B_h=0$ corrected solution only with \meqref{local1}.  }  \\
      \hline
 2&    0.93913487&0.0&      1.57251299&0.0&  1920\\
 3&    0.32826239&1.5&      1.03345712&0.6&  15360\\
 4&    0.05715300&2.5&      0.32677858&1.7&  122880\\
 5&    0.00665786&3.1&      0.02944335&3.5&983040\\
      \hline
\end{tabular}
 \caption{\lab{b-3} The errors and the order of convergence of $P_3$ solutions, for \meqref{sol1} with global $P_3 - H(div)$ and local interior modes only corrections.}
\end{center}

\end{table}

\section{Conclusion}

In this paper, we have proposed an efficient correction procedure to ensure
the divergence free condition of the magnetic field. The correction is done in
two steps: the first step can be done locally on each element which removes the
high order terms in the divergence error of the B-field; the second step is a
global one which removes the remaining constant term in the divergence error on
each element. Numerical results have shown the effectiveness of the proposed
method in enforcing divergence-free condition for a magnetic induction
equation while maintaining the accuracy of the solution itself.

\section*{Acknowledgement}

The author (W.C) acknowledges the support of the US Army Office of Research
(Grant No. W911NF-14-1-0297) and US National Science Foundation (Grant No.
DMS-1315128) and the National Natural Science Foundation of China (No.
91330110) for the work in this paper. The research of the author (J.H.) is
 supported by NSFC projects 11271035,  91430213 and 11421101.
Authors also like to thank Prof. Lin-bo Zhang for helpful discussions on implementation of H(div) basis.

\end{document}